\documentclass[12pt]{amsart}
\usepackage{amsmath,amsthm,amssymb}
\usepackage[arrow,matrix]{xy}
\usepackage{a4wide}
\usepackage{hyperref}
\usepackage{myheader}

\newcommand{\ZZ}{\mathbb{Z}_2}

\newcommand{\mathscr}{}

\begin{document}
\title {The number of small covers over cubes}

\author{Suyoung Choi}
\address{Department of Mathematical Science, Korea Advanced Institute of Science and Technology, Gu-sung Dong, Yu-sung Gu, Daejeon 305-701, Korea}
\thanks{This work was supported by the \textbf{SRC} Program of Korea Science and Engineering Foundation (KOSEF) grant funded by the Korea government(MOST) (No. R11-2007-035-02002-0).}

\email{choisy@kaist.ac.kr}

\date{\today}
\maketitle

\begin{abstract}
In the present paper we find a bijection between the set of small covers over an $n$-cube and the set of acyclic digraphs with $n$ labeled nodes. Using this, we give formulas of the number of small covers over an $n$-cube (generally, a product of simplices) up to Davis-Januszkiewicz equivalence classes and $\ZZ^n$-equivariant homeomorphism classes. Moreover we prove that the number of acyclic digraphs with $n$ unlabeled nodes is an upper bound of the number of small covers over an $n$-cube up to homeomorphism.
\end{abstract}

\tableofcontents

\section{Introduction}
Let $P$ be a simple convex polytope of dimension $n$ and $\mathcal{F}(P)=\{F_1 , \ldots, F_m\}$ be the set of facets of $P$.
Consider $\lambda:\mathcal{F}(P) \rightarrow \ZZ^n$ which satisfies the \emph{non-singularity condition}; $\{ \lambda(F_{i_1}) , \ldots , \lambda(F_{i_n}) \}$ is a basis of $\ZZ^n$ whenever the intersection $F_{i_1} \cap \cdots \cap F_{i_n}$ is non-empty. We call $\lambda$ a \emph{characteristic function}. Let $\ZZ(F_i)$ be the subgroup of $\ZZ^n$ generated by $\lambda(F_i)$.

Given a point $p \in P$, we denote by $G(p)$ the minimal face containing $p$ in its relative interior. Assume $G(p) = F_{j_1} \cap \cdots \cap F_{j_k}$. Then $\ZZ(G(p)) = \oplus_{i=1}^{k} \ZZ(F_{j_i})$. Note that $\ZZ(G(p))$ is a $k$-dimensional subgroup of $\ZZ^n$. Let $M(\lambda)$ denote $P \times (\ZZ)^n /\sim$, where $(p,g) \sim (q,h)$ if $p=q$ and $g^{-1}h \in \ZZ(G(p))$.
The free action of $\ZZ^n$ on $P \times \ZZ^n$  descends to an action on $M(\lambda)$ with quotient $P$. It is easy to show that the action is locally standard. On the other hand, Davis and Januszkiewicz introduced the notion of what is called a \emph{small cover} in \cite{DJ}. A small cover is a smooth closed manifold $M^n$ with a locally standard $\ZZ^n$-action such that its orbit space is a simple convex polytope. Thus, $M(\lambda)$ is a small cover over $P$.

Two small covers $M_1$ and $M_2$ are said to be \emph{weakly $\ZZ^n$-equivariantly homeomorphic} (or simply \emph{weakly $\ZZ^n$-homeomorphic}) if there is an automorphism $\varphi : \ZZ^n \rightarrow \ZZ^n$ and a homeomorphism $f:M_1 \rightarrow M_2$ such that $f(t \cdot x) = \varphi(t) \cdot f(x)$ for every $t \in \ZZ^n$ and $x \in M_1$. If $\varphi$ is an identity, then $M_1$ and $M_2$ are $\ZZ^n$-homeomorphic. Following Davis and Januszkiewicz, two small covers $M_1$ and $M_2$ over $P$ are said to be \emph{Davis-Januszkiewicz equivalent} (or simply, \emph{D-J equivalent}) if there is a weakly $\ZZ^n$-homeomorphism $f : M_1 \rightarrow M_2$ covering the identity on $P$.

Let $cf(P)$ denote the set of all characteristic functions over $P$.
\begin{theorem}[\cite{DJ}]
    All small covers over $P$ are given by $\{M(\lambda) | \lambda \in cf(P)\}$, i.e. for each small cover $M$ over $P$, there is a characteristic function $\lambda$ with $\ZZ^n$-homeomorphism $M(\lambda) \rightarrow M$ covering the identity on $P$.
\end{theorem}

There is a natural free left action of $GL(n, \ZZ)$ on $cf(P)$ defined by the correspondence $\lambda \mapsto \sigma \circ \lambda$. This action induces D-J equivalence on $cf(P)$. Hence the number of D-J equivalence classes over $P$ is $| GL(n, \ZZ) \setminus cf(P) |$. In recent years, numerous studies have attempted to enumerate the number of equivalence classes of all small covers over a specific polytope. In \cite{MR1937435}, Garrison and Scott used a computer program to show that the number of D-J classes over a dodecahedron is 2165. Moreover they calculated that the number of homeomorphism classes of all small covers over a dodecahedron is 25. In \cite{MR2328006}, Cai, Chen and L\"{u} calculated the number of D-J classes and $\ZZ^3$-equivariant homeomorphism classes over 3-dimensional prisms.

The arrangement of this paper is as follows. In \autoref{section:small covers over cubes} we study about small covers over cubes and some algebraic lemmas and use a combinatorial method to prove that the number of small covers over an $n$-cube up to D-J equivalence is equal to the number of acyclic digraphs with $n$ labeled nodes. Moreover we give the formula of the number of small covers over a product of simplices up to D-J equivalence. In \autoref{section:Equivariant homeomorphism classes}, we obtain a calculation formula of the number of equivariant homeomorphism classes of all small covers over cubes. In \autoref{section:upper bound}, we show that the number of weakly homeomorphism classes is less than or equal to the number of acyclic digraphs with $n$ unlabeled nodes.

\section{Small covers over cubes and acyclic digraphs} \label{section:small covers over cubes}
\subsection{Small covers over an $n$-cube}
Recall that we may assign an $(n \times m)$-matrix $\Lambda$ to an element $\lambda \in cf(P)$ by ordering the facets and choosing a basis for $(\ZZ)^n$.
$$
    \Lambda = \left( \lambda(F_1) \cdots \lambda(F_m) \right) = ( A | B ),
$$ where $A$ be an $(n \times n)$-matrix and $B$ be an $(n \times (m-n))$-matrix. Since there is a 1-1 correspondence between the D-J classes over $P$ and $GL(n,\ZZ)\setminus cf(P)$, up to D-J equivalence, the refined representative in its coset class is given by $\Lambda \sim (E_n | A^{-1}B)$, where $E_n$ is an identity matrix of size $n$. Denote $\Lambda_\ast=A^{-1}B$. We refer to $\Lambda$ as the \emph{refined form} of characteristic function $\lambda$, and call $\Lambda_\ast$ its \emph{reduced submatrix}.

When $P$ is an $n$-dimensional cube $I^n$, since the number of
facets of $I^n$ is $2n$, $\Lambda$ is an ($n \times 2n$)-matrix,
i.e., $\Lambda_\ast$ is an ($n\times n$)-matrix. We shall
additionally assume that the facets $F_j$ and $F_{n+j}$ do not
intersect for $1 \leq j \leq n$. Then the non-singularity
condition of $\lambda$ is equivalent to the following; \emph{every
principal minor of $\Lambda_\ast$ is 1}. Let $M(n)$ be the set of
$\ZZ$-matrices of size $n$ all of whose principal minors are 1.
Then we have the following 1-1 correspondence:
$$
    \{\text{D-J classes over }I^n\} \leftrightarrow M(n).
$$

Given a permutation $\mu$ of $n$ elements, denote by $P(\mu)$ the corresponding $n \times n$ permutation matrix, which has units in positions $(\mu(i),i)$ for $1 \leq i \leq n$, and zeros otherwise. There is the symmetric group action on $n \times n$ matrices by conjugations $A \mapsto P(\mu)^{-1} A P(\mu)$. Now we give the technical lemma which is first appeared in \cite{MR1850429}. But we cited it from \cite{MP}, Lemma 3.3.

\begin{lemma}[\cite{MP}] \label{lemma:technical lemma}
    Let $R$ be a commutative integral domain with an identity element 1, and let $A$ be an $n \times n$ matrix with entries in $R$. Suppose that every proper principal minor of $A$ is 1. If $\det A =1$, then $A$ is conjugate to a unipotent upper triangular matrix  by a permutation matrix, and otherwise to a matrix of the form
    $$
\left(
  \begin{array}{ccccc}
    1 & b_1 & 0 & \cdots & 0 \\
    0 & 1 & b_2 & \cdots & 0 \\
    \vdots & \vdots & \ddots & \ddots & \vdots \\
    0 & 0 & \cdots & 1 & b_{n-1} \\
    b_n & 0 & \cdots & 0 & 1 \\
  \end{array}
\right)
    $$
    where $b_i \neq 0$ for every $i$.
\end{lemma}

    \begin{theorem}\label{theorem:main theorem}
    The number of acyclic digraphs with $n$ labeled nodes is equal to the number of D-J equivalence classes of all small covers over $I^n$.
    \end{theorem}
    A ``digraph" means a graph with at most one edge directed from vertex $i$ to vertex $j$, for $1 \leq i \leq n, 1 \leq j \leq n$. An ``acyclic" means there is no cycle of any length.
    \begin{proof}
    Let $G$ be a digraph with $n$ labeled nodes. Let $A(G)$ be the vertex adjacency matrix of $G$ with $\ZZ$ entries.
    Set $B(G) := E_n + A(G)$. Note that conjugation action at $B(G)$ can be regarded as relabelling nodes.
    It is obvious that $G$ is acyclic if and only if $A(G)$ would be a strictly upper triangular matrix by conjugation,
    i.e. $B(G)$ would be an upper triangular matrix with diagonal entries 1.
    Note that the determinant of $B(G)$ is 1. Let $\mathcal{G}_n$ be the set of acyclic digraphs with labeled $n$ nodes.
    Define $\phi : \mathcal{G}_n \rightarrow M(n)$ by $G \mapsto B(G)$. We claim that it is indeed the bijection between $\mathcal{G}_n$ and $M(n)$. Let $G$ be an acyclic digraph with $n$ labeled nodes and $V$ be a node set of $G$. For any subset $V'$ of $V$, consider the induced subgraph $G'(V')$ by $V'$. Then it is also an acyclic digraph. If $|V'|=k$, then $B(G'(V'))$ is a $k$-rowed submatrix of $B(G)$. (A \emph{$k$-rowed submatrix} of $n \times n$ matrix $B$ is a $k \times k$ submatrix of $B$ whose entries, $b_{i,j}$, have indices $i$ and $j$ that are the elements of the same $k$-element subset of $\{1, \ldots, n\}$.) Thus the determinant of $B(G')$ is the principal minor of $B(G)$. Since $\det B(G') =1$ for any acyclic digraph $G'$, hence $B(G) \in M(n)$. Thus $\phi$ is well-defined.
    And it is injective since $G \mapsto A(G)$ is injective. Let $B$ be an element
of $M(n)$. Note that every diagonal entry of $B$ is 1. Then there
is a digraph $G$ such that $B(G)=B$. By \autoref{lemma:technical
lemma}, since every principal minor of $B$ is 1, $B$ is conjugated
to a unipotent upper triangular matrix. This implies $G$ is
acyclic, thus $\phi$ is surjective.
    \end{proof}

    Acyclic digraphs were counted by Robinson in \cite{MR0276143} and by Stanley in \cite{MR0317988}.

    \begin{theorem}[\cite{MR0276143}, \cite{MR0317988}]
    Let $R_n$ be the number of acyclic digraphs with $n$ labeled nodes. Then
    $$
        R_n = \sum_{k=1}^n (-1)^{k+1} {n \choose k} 2^{k(n-k)} R_{n-k}.
    $$
    \end{theorem}
$$
\begin{tabular}{|c|c c c c c c c c c}
  \hline
  $n$ & 0 & 1 & 2 & 3 & 4 & 5 & 6 & 7 & $\cdots$\\ \hline
  $R_n$ & 1 & 1 & 3 & 25 & 543 & 29281 & 3781503 & 1138779265 & $\cdots$ \\
  \hline
\end{tabular}
$$

Now consider $A \in M(n)$. Note that $A$ is a $\ZZ$-matrix. One may
regard $A$ as a real matrix with entries in $\{0, 1\}$. We simply
call it \emph{$(0,1)$-matrix}. Then every principal minor of $A$
is an odd number.
\begin{lemma} \label{lemma:principal minor = 1}
    Let $A$ be a $(0,1)$-matrix all of whose principal minors are odd. Then every principal minor of $A$ is 1.
\end{lemma}
\begin{proof}
We shall use an induction on $n$. When $n=1$, it is obvious. Assume that it holds for matrices of size $\leq n-1$. By induction hypothesis, every proper principal minor of $A$ is +1. If det$(A) \neq 1$, by \autoref{lemma:technical lemma}, $\det(A) = 1 \pm \prod b_i$. However $A$ is a (0,1)-matrix and a conjugation action is a permutation of the rows and columns of $A$, hence $b_i$'s must be 1. Thus, the determinant of $A$ is even. This is a contradiction. Thus $\det(A) = 1$.
\end{proof}

\begin{corollary}
    The number of acyclic digraphs with $n$ labeled nodes is equal to the number of real (0,1)-matrices all of whose principal minors are 1.
\end{corollary}
\begin{proof}
    This follows immediately from \autoref{lemma:principal minor = 1} and \autoref{theorem:main theorem}.
\end{proof}

\begin{remark}
We know that the number of acyclic digraphs with $n$ labeled nodes is equal to the number of (0,1)-matrices whose eigenvalues are positive. It is conjectured by Weisstein in 2001 and proved by McKay, \emph{et al}, in \cite{MR2085343}. Thus it can be easily checked that all eigenvalues of (0,1)-matrix are positive if and only if all of its principal minors are 1.
\end{remark}

\begin{remark}
We can define a \emph{quasitoric manifold} with $(S^1)^n$-action
as we did above. In this case, the reduced submatrix
$\Lambda_\ast$ of a characteristic function is an integer matrix.
When $P$ is an $n$-cube, $\Lambda_\ast$ is an ($n \times
n$)-matrix all of whose principal minors are $\pm 1$. Especially,
if every principal minor of $-\Lambda_\ast$ is 1, then the
quasitoric manifold is equivalent to a Bott tower. We refer the
reader \cite{MP}. By \autoref{lemma:principal minor = 1}, for a small cover over $I^n$ with $\Lambda_\ast$, there is the Bott tower over $I^n$ such that the characteristic function whose reduced submatrix is $-\Lambda_\ast$ as a (0,1)-matrix.
\end{remark}

\subsection{Small covers over a product of simplices}
The above processes can be extended to the case of a product of
simplices. Let $P = \prod_{i=1}^l \Delta^{n_i}$ with $\sum_{i=1}^l
n_i = n$, where $\Delta^{n_i}$ is the $n_i$-simplex for $i = 1,
\ldots, l$. Let $\{f_0^i , \ldots, f_{n_i}^i \}$ be the set of facets of the simplex $\Delta^{n_i}$. Therefore the set of facets of $P$ is
$$
    \{ F^i_{k_i} | 0 \leq k_i \leq n_i , i=1, \ldots, l \}
$$ where $F^i_{k_i} = \Delta^{n_1} \times \cdots \times \Delta^{n_{i-1}} \times f_{k_i}^i \times \Delta^{n_{i+1}} \times \cdots \times \Delta^{n_l}$. Thus there are $n+l$ facets in $P$.
 Then a reduced submatrix $\Lambda_\ast$ of characteristic function over $P$ is an ($n \times l$)-matrix. On the other hand, $\Lambda_\ast$ can be viewed as an ($l \times l$)-matrix ($\mathbf{v}_{i,j}$) whose entries in the $j$-th row are vectors in $\ZZ^{n_j}$. We shall call it a \emph{vector matrix}. We refer the reader \cite{ChoiMasudaSuh} for details. Let $\Lambda_{k_1 \cdots k_l}$ be the $(l \times l)$-submatrix of $\Lambda$ whose $j$-th row is the $k_j$-th row of the $\mathbf{v_{i,j}}$. Then the non-singularity condition for $\Lambda$ is equivalent that every principal minor of $\Lambda_{k_1 \cdots k_l}$ is 1 for any $1 \leq k_1 \leq n_1, \ldots, 1 \leq k_l \leq n_l$.

\begin{theorem}
Let $\sharp DJ( \prod_{i=1}^l \Delta^{n_i})$ denote the number of D-J equivalence classes over $\prod_{i=1}^l \Delta^{n_i}$. Then
$$\sharp DJ( \prod_{i=1}^l \Delta^{n_i}) = \sum_{G \in \mathcal{G}_l} \prod_{v_i \in V(G)} (2^{n_i}-1)^{\text{outdeg}(v_i)},$$
where $\mathcal{G}_l$ is the set of acyclic digraphs with labeled
$l$ nodes and $V(G)=\{v_1, \ldots, v_l\}$ is the labeled vertex
set of $G$.
\end{theorem}
\begin{proof}
Let $\Lambda_\ast = (\mathbf{v}_{i,j}$) be a reduced submatrix of characteristic function over $P$ with $\mathbf{v}_{i,j} \in \ZZ^{n_i}$.
Denote by $B(\Lambda_\ast):=(b_{i,j})$ the corresponding ($l\times l$)-matrix over $\ZZ$, which has units in positions $(i,j)$ if
    $\mathbf{v}_{i,j}$ is nonzero, and zeros otherwise. Define the map $\psi$ from $\{ GL(n,\ZZ) \setminus cf(P) \}$ to $\mathcal{G}_l$ by $\Lambda_\ast \mapsto G$ such that the adjacency matrix of $G$ is $B(\Lambda_\ast) - E_l$.
    By using similar arguments of the proof of \autoref{theorem:main
theorem}, we can prove that $\psi$ is well-defined.
Thus the number of characteristic functions is $\sum_{G \in
\mathcal{G}_l} | \psi^{-1}(G)|$. Let $G$ be an element of
$\mathcal{G}_l$ and $\Lambda_\ast = (\mathbf{v}_{i,j})$ be an $(l
\times l)$-vector matrix such with $\psi(\Lambda_\ast)=G$. Note
that an directed edge from $i$ to $j$ in $G$ is associated to a
nonzero $\mathbf{v}_{i,j}$ and $\mathbf{v}_{i,j} \in \ZZ^{n_i}$. Note that
$B(\Lambda_\ast)$ is conjugated to a unipotent upper triangular
matrix. Therefore the non-singularity condition holds for
arbitrary nonzero vectors in non-diagonal entries, i.e. we have $2^{n_i}-1$ choices for each nonzero vector $\mathbf{v}_{i,j}$. Thus
$|\psi^{-1}(G)| = \prod_{e \in E(G)} (2^{n_{i(e)}}-1) =
\prod_{v_i \in V(G)}(2^{n_i}-1)^{\text{outdeg}(v_i)}$, where $E(G)$ is the set of directed edges of $G$ and $i(e)$ is the index of the initial vertex of $e \in E(G)$.
\end{proof}

\begin{example}
\underline{$l =2$} : $ \sharp DJ( \Delta^{n_1} \times \Delta^{n_2}) = 1 + (2^{n_1} -1) +(2^{n_2} - 1)$.

\underline{$l =3$} : $ \sharp DJ( \Delta^{n_1} \times \Delta^{n_2} \times \Delta^{n_3}) = 1 + 2(x_1 + x_2 + x_3) + (x_1+x_2+x_3)^2+ (x_1 x_2 + x_2 x_3 + x_3 x_1) +  (x_1 + x_2 + x_3)(x_1^2 + x_2 ^2 + x_3 ^2) - x_1^3- x_2^3- x_3^3$, where $x_i = 2^{n_i}-1$ for $i=1,2,3$.
\end{example}

\section{Counting $\ZZ^n$-equivariant homeomorphism classes}\label{section:Equivariant homeomorphism classes}
Let $P$ be a simple convex polytope of dimension $n$ and $\mathcal{F}(P)$ be the set of faces of $P$. An \emph{automorphism} of $\mathcal{F}(P)$ is a bijection from $\mathcal{F}(P)$ to itself which preserves the poset structure of all faces of $P$. Let $\text{Aut}(\mathcal{F}(P))$ denote the group of automorphisms of $\mathcal{F}(P)$. One can define the right action of $\text{Aut}(\mathcal{F}(P))$ on $cf(P)$ by $\lambda \times h \mapsto \lambda \circ h$, where $\lambda \in cf(P)$ and $h \in \text{Aut}(\mathcal{F}(P))$. The following theorem is well-known. We refer the reader \cite{LM}.
\begin{theorem}\label{theorem:equivariant homeomorphism}
    Two small covers over an $n$-dimensional simple convex polytope $P$ are $\ZZ^n$-equivariantly homeomorphic if and only if there is $h \in \text{Aut}(\mathcal{F}(P))$ such that $\lambda_1 = \lambda_2 \circ h$, where $\lambda_1$ and $\lambda_2$ are characteristic functions of small covers.
\end{theorem}

Thus we are going to count the orbits of $cf(I^n)$  under the action of $Aut(\mathcal{F}(I^n))$. The Burnside's formula is very useful in the enumeration of the number of orbits.

\begin{lemma} [Burnside's formula]
    Let $G$ be a finite group acting on a set $X$. Then the number of orbits of $X$ under the $G$-action is equal to $\frac{1}{|G|} \sum_{g \in G} |X^g|,
    $ where $X^g = \{x \in X | gx=x\}$.
\end{lemma}

\begin{theorem}
Let $Q_n$ be the number of $\ZZ^n$-equivariant homeomorphism classes of small covers over $I^n$ and $R_k$ be the number of acyclic digraphs with $k$ labeled nodes.Then
    $$ Q_n = \frac{\sum_{k=0}^n {n \choose k}2^{k(n-k)}R_k}{2^n n!} \cdot \prod_{i=0}^{n-1}(2^n-2^i) $$
\end{theorem}

\begin{proof}
All elements of $\text{Aut}(\mathcal{F}(I^n))$ can be written in a simple form as follows :
$$
    \mu \cdot \chi_1^{e_1} \cdot \cdots \cdot \chi_n^{e_n} , e_j \in \ZZ
$$ with a permutation $\mu \in S_n$ and reflections $\chi_1, \ldots, \chi_n$. Hence $|\text{Aut}(\mathcal{F}(I^n))| = 2^n n!$. Note that $\mu$ is a permutation of the pairs of opposite facets and $\chi_i$ is the interchange of $i$-th opposite facets for each $i$. For some $g=\mu \chi_1^{e_1} \cdots \chi_n^{e_n} \in \text{Aut}(\mathcal{F}(I^n))$, let $cf(I^n)^g$ denote the set of elements in $cf(I^n)$ fixed by $g$. First, we claim that $cf(I^n)^g$ is nonempty implies $\mu = 1$. Let $\lambda \in \text{Aut}(\mathcal{F}(I^n))^g$ and $\mathcal{F}(I^n) = \{ F_1, \ldots, F_{2n} \}$ be the set of facets of $I^n$ such that $F_i \cap F_{n+i} = \emptyset$ for all $i=1, \ldots, n$. Note that the non-singularity condition implies the determinant of  $(\lambda(F_{\epsilon(1)}) \cdots \lambda(F_{\epsilon(n)}))$ is 1, where $\epsilon(t)$ is either $t$ or $n+t$. Thus there is no pair $F_i,F_j$ such that $\lambda(F_i) = \lambda(F_j)$ and $n \nmid i-j$. We deduce that $\mu=1$. Now, we are going to enumerate $|cf(I^n)^g|$ when $\mu =1$. We may assume $g=\chi_{1} \cdots \chi_{k}$ for some $k$. Let $\lambda$ be an element of $cf(I^n)^g$ and $\Lambda$ be an ($n\times 2n$)-matrix corresponding to $\lambda$. Note that $\Lambda = (A|B) = A \cdot (E_n | \Lambda_\ast)$, where $\Lambda_\ast = A^{-1}B$. Note that $\lambda$ is fixed by $g$ if and only if the first $k$ columns of $A$ and $B$ are the same. Thus $\Lambda_\ast$ is of the following form:
$$
    \left(
         \begin{array}{cc}
            E_k & S \\
            0 & T \\
         \end{array}
    \right),
$$ where $E_k$ is an identity matrix of size $k$, $T$ is an $((n-k)\times (n-k))$-matrix and $S$ is a $(k \times (n-k))$-matrix. Note that $\Lambda_\ast \in M(n)$ if and only if $T \in M(k)$. This implies $|cf(I^n)^g|= |GL(n,\ZZ)| \times 2^{k(n-k)} R_k$. Note that $|cf(I^n)^g|$ is independent of choices of $k$ $\chi_i$'s.
Thus, by Burnside's formula,
$$
   Q_n = \frac{\sum_{k=0}^n {n \choose k}2^{k(n-k)} R_k}{|\text{Aut}(I^n)|} \cdot |GL(n,\ZZ)|.
$$
The theorem is proved with well-known fact $|GL(n,\ZZ)| = \prod_{i=0}^{n-1}(2^n -2^i)$.
\end{proof}
$$
\begin{tabular}{|c|c c c c c c c }
  \hline
  $n$ & 0& 1& 2 & 3 & 4 & 5 & $\cdots$\\ \hline
  $Q_n$ & 1 & 1 & 6 & 259 & 87360 & 236240088 & $\cdots$ \\
  \hline
\end{tabular}
$$

\section{Upper bounds of the numbers of homeomorphism classes}\label{section:upper bound}
By \autoref{theorem:equivariant homeomorphism}, we have a one-to-one correspondence between the set of weakly equivariant homeomorphism classes of small covers over simple polytope $P$ and the double coset class by $GL(n,\ZZ)$ and $\text{Aut}(\mathcal{F}(P))$ on $cf(P)$. Let $T_n$ be the number of weakly $\ZZ^n$-equivariant homeomorphism classes of small covers over $I^n$. Then,
$$
    T_n = | GL(n,\ZZ) \setminus cf(I^n) / \text{Aut}(\mathcal{F}(I^n)) | = | M(n) / \text{Aut}(\mathcal{F}(I^n)) |
$$where $M(n)$ is the set of $\ZZ$-matrices of size $n$ all of whose principal minors are 1. Recall that $\text{Aut}(\mathcal{F}(I^n))$ consists of elements of the form $\mu \cdot \chi_1^{e_1} \cdot \cdots \cdot \chi_n^{e_n} , e_j \in \ZZ $ with a permutation $\mu \in S_n$ and reflections $\chi_1, = \cdots = \chi_n$. Consider the permutation group $S_n = \{ g \in \text{Aut}(\mathcal{F}(I^n)) | g=\mu \chi^0 \cdots \chi^0 \}$ as a subgroup of $\text{Aut}(\mathcal{F}(I^n))$. Then the action of $S_n$ on the set of facets of $I^n$ by permuting the pairs of opposite facets. Let $\Lambda$ be an $(n \times 2n)$ characteristic matrix. Then $\mu \in S_n$ acts as
$$
    \Lambda \mapsto \Lambda \cdot \left(
                                   \begin{array}{cc}
                                     P(\mu) & 0 \\
                                     0 & P(\mu) \\
                                   \end{array}
                                 \right).
$$

Thus, $(E_n | \Lambda_\ast) \mapsto (P(\mu) | \Lambda_\ast P(\mu)) \sim (E_n | P(\mu)^{-1} \Lambda_\ast P(\mu))$. This implies the action of $S_n$ on $M(n)$ is the conjugation action. That is a relabeling on nodes of acyclic digraphs. Hence we have the following theorem and corollary:

\begin{theorem}
    The number of weakly $\ZZ^n$-equivariant homeomorphism classes of small covers over $I^n$ is less than or equal to the number of acyclic digraphs with $n$ unlabeled nodes.
\end{theorem}

\begin{corollary}
    The number of homeomorphism classes of small covers over $I^n$ is less than or equal to the number of acyclic digraphs with $n$ unlabeled nodes.
\end{corollary}

Acyclic digraphs with unlabeled nodes were counted by Robinson in \cite{MR0476569}.
$$
\begin{tabular}{|c|c c c c c c c c c}
  \hline
  $n$ & 0 & 1 & 2 & 3 & 4 & 5 & 6 & 7 & $\cdots$\\ \hline
  $T_n$ & 1 & 1 & 2 & 6 & 31 & 302 & 5984 & 243668 & $\cdots$ \\
  \hline
\end{tabular}
$$

\section*{Acknowledgements}
I am grateful to my advisor, Dong Youp Suh, for his encouragement and a number of comments and suggestions. And I thank Jang Soo Kim for his comments, which improved the proof of \autoref{theorem:main theorem}. I also thank Mikiya Masuda for his suggestions, which pointed out about the difficulty in the smooth category for small covers.

\bigskip
\bibliographystyle{amsalpha}
\bibliography{reference}

\end{document}